\documentclass[11pt,twoside]{article}

\usepackage[latin1]{inputenc}
\usepackage{epsfig}
\usepackage{amsfonts}
\usepackage{cite}
\usepackage{color}
\definecolor{lgrey}{gray}{0.99}

\title{
  {\huge Complex Analysis of Real Functions \\[1.5ex]}
  I: Complex-Analytic Structure and \\
  Integrable Real Functions }

\author{
  \Large Jorge L. deLyra\footnote{Email: delyra@latt.if.usp.br} \\
  Department of Mathematical Physics \\
  Physics Institute \\
  University of São Paulo }

\date{May 26, 2018}

\newcommand{\ii}{\mbox{\boldmath$\imath$}}
\newcommand{\iii}{\mbox{\boldmath\scriptsize$\imath$}}

\newcommand{\e}[1]{\,{\rm e}^{#1}}

\newcommand{\ldot}{\mbox{\Large$\cdot$}\!}

\newtheorem{definition}{Definition}
\newtheorem{property}{Property}[definition]

\newtheorem{theorem}{Theorem}

\newtheorem{proof}{Proof}[theorem]
\newcommand{\Colon}{{\hspace{-0.4em}\bf:}}

\begin{document}\maketitle

\begin{abstract}
  \noindent
  A complex-analytic structure within the unit disk of the complex plane
  is presented. It can be used to represent and analyze a large class of
  real functions. It is shown that any integrable real function can be
  obtained by means of the restriction of an analytic function to the unit
  circle, including functions which are non-differentiable, discontinuous
  or unbounded. An explicit construction of the analytic functions from
  the corresponding real functions is given. The complex-analytic
  structure can be understood as an universal regulator for analytic
  operations on real functions.
\end{abstract}

\section{Introduction}\label{Sec01}

In this paper we will exhibit a mathematical structure, based on certain
analytic functions within the unit circle of the complex plane, that can
be used to represent and analyze a very wide class of real functions.
These include analytic and non-analytic integrable real functions, as well
as unbounded integrable real functions. All these objects will be
interpreted as parts of a larger complex-analytic structure, within which
they can be treated and manipulated in a robust and unified way.

In order to assemble the mathematical structure a set of mathematical
objects must be introduced, and their properties established. This will be
done in Section~\ref{Sec02}, in which all the eight necessary definitions
will be given, and all the corresponding properties will be stated and
proved. The objects to be defined are elements within complex
analysis~\cite{CVchurchill}, and include a general scheme for the
classification of all possible singularities of analytic functions, as
well as the concept of infinite integral-differential chains of functions.

As a first and important application of this complex-analytic structure,
in Section~\ref{Sec03} we will establish the relation between the
complex-analytic structure and integrable real functions. There we will
show that every integrable real function defined within a finite interval
corresponds to an inner analytic function and can be obtained by means of
the restriction of the real part of that analytic function to the unit
circle of the complex plane.

This is the first of a series of papers. The discussion of some parts and
aspects of this line of work will be postponed to forthcoming papers, in
order to keep each paper within a reasonable length. In the second paper
of the series we will extend the complex-analytic structure presented in
this paper, to include the whole space of singular Schwartz distributions,
also known as generalized real functions.

In the third paper of the series we will show that the whole Fourier
theory of integrable real functions is contained within that same
complex-analytic structure. We will show that this structure induces a
very general and powerful summation rule for Fourier series, that can be
used to add up Fourier series in a consistent way, even when they are
explicitly and strongly divergent. The complex-analytic structure will
then allow us to extend the Fourier theory beyond the realm of integrable
real functions, with the use of that summation rule.

In the fourth paper of the series we will show that one can include in the
same complex-analytic structure a large class of non-integrable real
functions, among those that are locally integrable almost everywhere. We
will see that the complex-analytic structure allows us to associate to
each such function a definite set of Fourier coefficients, despite the
fact that the functions are not integrable on the unit circle. There are
also other applications of the structure discussed here, for example in
the two-dimensional Dirichlet problem in partial differential equations, a
discussion of which will be given in the fifth paper of the series.

The material contained in this paper is a development, reorganization and
extension of some of the material found, sometimes still in rather
rudimentary form, in the
papers~\cite{FTotCPI,FTotCPII,FTotCPIII,FTotCPIV,FTotCPV}.

\section{Definitions and Properties}\label{Sec02}

Here we will introduce the definitions and basic properties of some
objects and structures which are not usually discussed in complex
analysis~\cite{CVchurchill}, and which we will use in the subsequent
sections. Consider then the unit circle centered at the origin of the
complex plane. Its interior is the open unit disk we will often refer to
along the paper. Any reference to the unit disk or to the unit circle
should always be understood to refer to those centered at the origin.

\begin{definition}\Colon\label{Def01}
  Inner Analytic Functions
\end{definition}

\noindent
A complex function $w(z)$ which is analytic in the open unit disk will be
named an {\em inner analytic function}. We will consider the set of all
such functions. We will also consider the subset of such functions that
have the additional property that $w(0)=0$, which we will name {\em proper
  inner analytic functions}.

\vspace{2.6ex}

\noindent
Note, in passing, that the set of all inner analytic functions forms a
vector space over the field of complex numbers, and so does the subset of
all proper inner analytic functions.

The focus of this study will be the set of real objects which are obtained
from the real parts of these inner analytic functions when we take the
limit from the open unit disk to its boundary, that is, to the unit
circle. Specifically, if we describe the complex plane with polar
coordinates $(\rho,\theta)$, then an arbitrary inner analytic function can
be written as

\begin{equation}
  w(z)
  =
  u(\rho,\theta)+\ii v(\rho,\theta),
\end{equation}

\noindent
where $\ii$ is the imaginary unit, and we consider the set of real objects
$f(\theta)$ obtained from the set of all inner analytic functions as the
limits of their real parts, from the open unit disk to the unit circle,

\begin{equation}
  f(\theta)
  =
  \lim_{\rho\to 1_{(-)}}u(\rho,\theta),
\end{equation}

\noindent
when and where such limits exist, or at least can be defined in a
consistent way.

Note that an inner analytic function may have any number of singularities
on the unit circle, as well as in the region outside the unit circle. The
concept of a singularity is the usual one in complex analysis, namely that
a singular point is simply a point where the function fails to be
analytic. The singularities on the unit circle will play a particularly
important role in the complex-analytic structure to be presented in this
paper. If any of these singularities turn out to be branch points, then we
assume that the corresponding branch cuts extend outward from the unit
circle, either out to infinity or connecting to some other singularity
that may exist outside the open unit disk.

Note also that the imaginary parts of the inner analytic functions do {\em
  not} generate an independent set of real objects, since the imaginary
part $v(\rho,\theta)$ of the inner analytic function $w(z)$ is also the
real part of the inner analytic function $\bar{w}(z)$ given by

\begin{equation}
  \bar{w}(z)
  =
  -\ii w(z).
\end{equation}

\noindent
We thus see, however, that the inner analytic functions do organize the
real functions in matched pairs, those originating from the real and
imaginary parts of each inner analytic function. The two real functions
forming such a pair may be described as mutually {\em Fourier conjugate}
functions. Finally, we will assume that, at all singular points where the
functions $w(z)$ can still be defined by continuity, they have been so
defined.

In addition to establishing this correspondence between complex functions
on the unit disk and real function of the unit circle, we will find it
necessary to define analytic operations on the complex functions that
correspond to the ordinary operations of differentiation and integration
on the real functions. As will be shown in what follows, the next two
definitions accomplish this.

\begin{definition}\Colon\label{Def02}
  Angular Differentiation
\end{definition}

\noindent
Given an arbitrary inner analytic function $w(z)$, its {\em angular
  derivative} is defined by

\begin{equation}
  w^{\ldot}(z)
  =
  \ii
  z\,
  \frac{dw(z)}{dz}.
\end{equation}

\noindent
The angular derivative of $w(z)$ will be denoted by the shifted dot, as
shown. The second angular derivative will be denoted by $w^{2\ldot}(z)$,
and so on.

\vspace{2.6ex}

\noindent
Note that this definition has been tailored in order for the following
property to hold.

\begin{property}\Colon\label{Prop1LogDif}
  In terms of the variables $(\rho,\theta)$, angular differentiation is
  equivalent to partial differentiation with respect to $\theta$, taken at
  constant $\rho$.
\end{property}

\noindent
Writing $z=\rho\exp(\ii\theta)$, and considering the partial derivative of
$z$ with respect to $\theta$, we have

\noindent
\begin{eqnarray}
  w^{\ldot}(z)
  & = &
  \ii\rho\e{\iii\theta}\,
  \frac{1}{\ii\rho\e{\iii\theta}}\,
  \frac{\partial w(z)}{\partial\theta}
  \nonumber\\
  & = &
  \frac{\partial u(\rho,\theta)}{\partial\theta}
  +
  \ii\,
  \frac{\partial v(\rho,\theta)}{\partial\theta},
\end{eqnarray}

\noindent
which establishes this property.

\vspace{2.6ex}

\noindent
Note that by construction we always have that $w^{\ldot}(0)=0$, so that we
may say that the operation of angular differentiation projects the space
of inner analytic functions onto the space of proper inner analytic
functions. We may now prove an important property of the angular
derivative.

\begin{property}\Colon\label{Prop2LogDif}
  The angular derivative of an inner analytic function is also an inner
  analytic function.
\end{property}

\noindent
Let us recall that the derivative of an analytic function always exists
and is also analytic, in the same domain of analyticity of the original
function. Since the constant function $w(z)\equiv\ii$ and the identity
function $w(z)\equiv z$ are analytic in the whole complex plane, and since
the product of analytic functions is also an analytic function, in their
common domain of analyticity, it follows at once that the angular
derivative of an inner analytic function is an inner analytic function as
well, which establishes this property.

\vspace{2.6ex}

\noindent
In other words, the operation of angular differentiation stays within the
space of inner analytic functions. Note that, since $w^{\ldot}(0)=0$,
angular differentiation always results in {\em proper} inner analytic
functions, and therefore that this operation also stays within the space
of proper inner analytic functions.

\begin{definition}\Colon\label{Def03}
  Angular Integration
\end{definition}

\noindent
Given an arbitrary inner analytic function $w(z)$, its {\em angular
  primitive} is defined by

\begin{equation}\label{EQLogInt}
  w^{-1\ldot}(z)
  =
  -\ii
  \int_{0}^{z}dz'\,
  \frac{w(z')-w(0)}{z'},
\end{equation}

\noindent
where the integral is taken along any simple curve from $0$ to $z$
contained within the open unit disk. Since the integrand is analytic
inside the open unit disk, including at the origin, as we will see shortly
while proving Property~\ref{Prop2LogInt}, due to the Cauchy-Goursat
theorem the integral does not depend on the curve along which it is taken.
The angular primitive will also be denoted by a shifted dot, this time
preceded by a negative integer, as indicated above.

\vspace{2.6ex}

\noindent
Let us prove that the apparent singularity of the integrand at $z=0$ is in
fact a removable singularity, so that the integrand can be defined at the
origin by continuity, thus producing a function which is continuous and
well defined there. If we simply take the $z\to 0$ limit of the integrand
we get

\begin{equation}
  \lim_{z\to 0}
  \frac{w(z)-w(0)}{z}
  =
  \frac{dw}{dz}(0),
\end{equation}

\noindent
since this limit is the very definition of the derivative of $w(z)$ at
$z=0$. Since $w(z)$ is an inner analytic function, and is thus analytic in
the open unit disk, it is differentiable at the origin, so that this limit
exists and is finite. Therefore the integrand can be defined at the origin
to have this particular value, so that it is continuous there. We assume
that the integrand is so defined at $z=0$, as part of the definition of
the angular primitive.

\vspace{2.6ex}

\noindent
Note that this definition has been tailored in order for the following
property to hold.

\begin{property}\Colon\label{Prop1LogInt}
  In terms of the variables $(\rho,\theta)$, angular integration can be
  understood as integration with respect to $\theta$, taken at constant
  $\rho$, up to an integration constant.
\end{property}

\noindent
Given any point $z$ in the open unit disk, and considering that we are
free to choose the path of integration from $0$ to $z$, we now choose to
go first from the origin along the positive real axis, until we reach the
radius $\rho$, and then to go along an arc of circle of radius $\rho$,
until we reach the angle $\theta$, thus separating the integral in two. In
the first integral the variations of $z$ are given by $dz=d\rho$, and in
the second one they are given by $dz=\ii\rho\exp(\ii\theta)d\theta$. Note
that as the integrand in Equation~(\ref{EQLogInt}) we have the {\em
  proper} inner analytic function given by

\noindent
\begin{eqnarray}
  w_{p}(z)
  & = &
  w(z)-w(0)
  \nonumber\\
  & = &
  u_{p}(\rho,\theta)
  +
  \ii
  v_{p}(\rho,\theta),
\end{eqnarray}

\noindent
where $w_{p}(0)=0$. The integral in Equation~(\ref{EQLogInt}) can
therefore be written as

\noindent
\begin{eqnarray}
  w^{-1\ldot}(z)
  & = &
  -\ii
  \int_{0}^{\rho}d\rho'\,
  \frac{w_{p}(\rho',0)}{\rho'}
  -
  \ii
  \int_{0}^{\theta}d\theta'\,
  \ii\rho\e{\iii\theta'}\,
  \frac{w_{p}(\rho,\theta')}{\rho\e{\iii\theta'}}
  \nonumber\\
  & = &
  C(\rho)
  +
  \int_{0}^{\theta}d\theta'\,
  u_{p}(\rho,\theta')
  +
  \ii
  \int_{0}^{\theta}d\theta'\,
  v_{p}(\rho,\theta'),
\end{eqnarray}

\noindent
where, in relation to the variable $\theta$, the integral on $\rho'$
becomes the complex function $C(\rho)$, which depends only on $\rho$ and
not on $\theta$, while the integral on $\theta'$ determines primitives
with respect to $\theta$ of the real and imaginary parts of $w_{p}(z)$,
which thus establishes this property.

\vspace{2.6ex}

\noindent
Note that by construction we always have that $w^{-1\ldot}(0)=0$, so that
we may say that the operation of angular integration projects the space of
inner analytic functions onto the space of proper inner analytic
functions. We may now establish an important property of the angular
primitive.

\begin{property}\Colon\label{Prop2LogInt}
  The angular primitive of an inner analytic function is also an inner
  analytic function.
\end{property}

\noindent
In order to prove that $w^{-1\ldot}(z)$ is an inner analytic function, we
use the power-series representation of the inner analytic function $w(z)$.
Since this function is analytic within the open unit disk, its Taylor
series around $z=0$, which is given by

\begin{equation}
  w(z)
  =
  w(0)
  +
  \sum_{k=1}^{\infty}
  c_{k}z^{k},
\end{equation}

\noindent
where $c_{k}=w^{k\prime}(0)/k!$ are the Taylor coefficients of $w(z)$ with
respect to the origin and where $w^{k\prime}(z)$ is the $k^{\rm th}$
derivative of $w(z)$, converges within that disk. We therefore have for
the integrand in Equation~(\ref{EQLogInt}) the power-series representation

\noindent
\begin{eqnarray}
  \frac{w(z)-w(0)}{z}
  & = &
  \sum_{k=1}^{\infty}
  c_{k}z^{k-1}
  \nonumber\\
  & = &
  \sum_{k=0}^{\infty}
  c_{k+1}z^{k}.
\end{eqnarray}

\noindent
Since this series has the same set of coefficients as the convergent
series of $w(z)$, it is equally convergent, as is implied for example by
the ratio test. Note that this shows, in particular, that the integrand is
analytic at $z=0$. Being a convergent power series, this series can be
integrated term by term, resulting in an equally convergent power series,
so that we have for the angular primitive

\noindent
\begin{eqnarray}
  w^{-1\ldot}(z)
  & = &
  -\ii
  \int_{0}^{z}dz'\,
  \sum_{k=0}^{\infty}
  c_{k+1}(z')^{k}
  \nonumber\\
  & = &
  -\ii
  \sum_{k=0}^{\infty}
  \frac{c_{k+1}}{k+1}\,
  z^{k+1}
  \nonumber\\
  & = &
  -\ii
  \sum_{k=1}^{\infty}
  \frac{c_{k}}{k}\,
  z^{k}.
\end{eqnarray}

\noindent
Due to the factors of $1/k$, when $k\to\infty$ the coefficients of this
series go to zero faster than those of the convergent Taylor series of
$w(z)$, and thus it is also convergent, in the same domain of convergence
of the Taylor series of $w(z)$. This confirms that this series is
convergent within the open unit disk. Being a convergent power series, it
converges to an analytic function, thus proving that $w^{-1\ldot}(z)$ is
analytic within the open unit disk. We may conclude therefore that the
angular primitive of an inner analytic function is an inner analytic
function as well, which establishes this property.

\vspace{2.6ex}

\noindent
In other words, the operation of angular integration stays within the
space of inner analytic functions. Note that, since $w^{-1\ldot}(0)=0$,
angular integration always results in {\em proper} inner analytic
functions, and therefore that this operation also stays within the space
of proper inner analytic functions.

\vspace{2.6ex}

\noindent
Let us now prove that the operations of angular differentiation and of
angular integration are inverse operations to one another. Strictly
speaking, this is true within the subset of inner analytic functions that
have the additional property that $w(0)=0$, that is, for proper inner
analytic functions. Since any inner analytic function can be obtained from
a proper inner analytic function by the mere addition of a constant, this
is a very weak limitation. Let us consider then the space of proper inner
analytic functions.

\begin{property}\Colon\label{Prop3LogInt}
  The angular primitive of the angular derivative of a proper inner
  analytic function is that same proper inner analytic function.
\end{property}

\noindent
We simply compose the two operations in the required order, and calculate
in a straightforward manner, merely using the fundamental theorem of the
calculus, to get

\noindent
\begin{eqnarray}
  -\ii
  \int_{0}^{z}dz'\,
  \frac{1}{z'}\,
  \left[
    \ii
    z'\,
    \frac{dw}{dz'}(z')
    -
    \ii
    0\,
    \frac{dw}{dz'}(0)
  \right]
  & = &
  \int_{0}^{z}dz'\,
  \frac{dw}{dz'}(z')
  \nonumber\\
  & = &
  w(z)-w(0),
\end{eqnarray}

\noindent
which is the original inner analytic function $w(z)$ so long as $w(0)=0$,
that is, for proper inner analytic functions, thus establishing this
property.

\begin{property}\Colon\label{Prop4LogInt}
  The angular derivative of the angular primitive of a proper inner
  analytic function is that same proper inner analytic function.
\end{property}

\noindent
We simply compose the two operations in the required order, and calculate
in a straightforward manner, using this time the power-series
representation of the inner analytic function $w(z)$. First integrating
term by term and then differentiating term by term, both of which are
allowed operations for convergent power series, we get

\noindent
\begin{eqnarray}
  \ii
  z\,
  \frac{d}{dz}
  (-\ii)
  \int_{0}^{z}dz'\,
  \frac{w(z')-w(0)}{z'}
  & = &
  z\,
  \frac{d}{dz}
  \int_{0}^{z}dz'\,
  \sum_{k=1}^{\infty}
  c_{k}(z')^{k-1}
  \nonumber\\
  & = &
  z\,
  \frac{d}{dz}
  \sum_{k=1}^{\infty}
  \frac{c_{k}}{k}\,z^{k}
  \nonumber\\
  & = &
  \sum_{k=1}^{\infty}
  c_{k}z^{k}
  \nonumber\\
  & = &
  w(z)-w(0),
\end{eqnarray}

\noindent
which is the original inner analytic function $w(z)$ so long as $w(0)=0$,
that is, for proper inner analytic functions, thus establishing this
property.

\vspace{2.6ex}

\noindent
With the use of the operations of angular differentiation and angular
integration the space of proper inner analytic functions can now be
organized as a collection of infinite discrete chains of functions, so
that within each chain the functions are related to each other by either
angular integrations or angular differentiations. This leads to the
definition that follows.

\begin{definition}\Colon\label{Def04}
  Integral-Differential Chains
\end{definition}

\noindent
Starting from an arbitrary proper inner analytic function $w(z)$, also
denoted as $w^{0\ldot}(z)$, one proceeds in the differentiation direction
to the functions $w^{1\ldot}(z)$, $w^{2\ldot}(z)$, $w^{3\ldot}(z)$, etc,
and in the integration direction to the functions $w^{-1\ldot}(z)$,
$w^{-2\ldot}(z)$, $w^{-3\ldot}(z)$, etc. One thus produces an infinite
chain of proper inner analytic functions such as

\begin{equation}
  \left\{
    \ldots,
    w^{-3\ldot}(z),
    w^{-2\ldot}(z),
    w^{-1\ldot}(z),
    w^{0\ldot}(z),
    w^{1\ldot}(z),
    w^{2\ldot}(z),
    w^{3\ldot}(z),
    \ldots\;
  \right\},
\end{equation}

\noindent
in which angular differentiation takes one to the right and angular
integration takes one to the left. We name such a structure an {\em
  integral-differential chain} of proper inner analytic functions. We may
refer to the proper inner analytic functions forming the chain as {\em
  links} in that chain.

\vspace{2.6ex}

\noindent
Note that all the functions in such a chain have exactly the same set of
singular points on the unit circle, although the character of these
singularities will change from function to function along the chain. Note
also that each such integral-differential chain induces, by means of the
real parts of their inner analytic functions, a corresponding chain of
real objects over the unit circle, when and where the limits from the open
unit disk to the unit circle exist, or can be consistently defined.
Finally note that, given a singularity at a certain point on the unit
circle, the integral-differential chain also induces a corresponding chain
of singularities at that point. Let us now prove an important property of
these chains, namely that they do not intersect each other.

\begin{property}\Colon\label{Prop1InDiCh}
  Two different integral-differential chains of proper inner analytic
  functions cannot have a member-function in common.
\end{property}

\noindent
In order to prove this, we start by proving that, if two proper inner
analytic functions have the same angular derivative, then they must be
equal. If we have two such proper inner analytic functions $w_{1}(z)$ and
$w_{2}(z)$, the statement that they have the same angular derivative is
expressed as

\noindent
\begin{eqnarray}
  \ii
  z\,
  \frac{d}{dz}w_{2}(z)
  -
  \ii
  z\,
  \frac{d}{dz}w_{1}(z)
  & = &
  0
  \;\;\;\Rightarrow
  \nonumber\\
  \frac{d}{dz}[w_{2}(z)-w_{1}(z)]
  & = &
  0
  \;\;\;\Rightarrow
  \nonumber\\
  w_{2}(z)-w_{1}(z)
  & = &
  C,
\end{eqnarray}

\noindent
where $C$ is some complex constant, for all $z$ within the open unit disk,
including the case $z=0$, as one can see if one takes the limit $z\to 0$
of the last equation above. However, since at $z=0$ we have that
$w_{1}(0)=0$ and $w_{2}(0)=0$, it then follows that $C=0$, so that we may
conclude that

\begin{equation}
  w_{2}(z)
  \equiv
  w_{1}(z),
\end{equation}

\noindent
thus proving the point. A similar result is valid for two proper inner
analytic functions that have the same angular primitive. Since we have
already shown that angular integration and angular differentiation are
inverse operations to each other, we can prove this by simply noting the
trivial fact that the operation of angular differentiation cannot produce
two different results for the same function. Therefore, there cannot exist
two different proper inner analytic functions whose angular primitives are
one and the same function.

We may now conclude that two different integral-differential chains of
proper inner analytic functions can never have a member-function in
common, because this would mean that two different proper inner analytic
functions would have either the same angular derivative or the same
angular primitive, neither of which is possible. It follows that each
proper inner analytic function appears in one and only one of these
integral-differential chains, which establishes this property.

\vspace{2.6ex}

\noindent
Note, for future use, that there is a {\em single} integral-differential
chain of proper inner analytic functions which is a constant chain, in the
sense that all member-functions of the chain are equal, namely the null
chain, in which all members are the null function $w(z)\equiv 0$. It is
easy to verify that the differential equation $w^{\ldot}(z)=w(z)$ has no
other inner analytic function as a solution. Note also that one may
consider all the non-proper inner analytic functions $w(z)$ which are
related to a given proper inner analytic function $w_{p}(z)$ to also
belong to the same link of the corresponding integral-differential
chain. Since all such functions have the form

\begin{equation}
  w(z)
  =
  C
  +
  w_{p}(z),
\end{equation}

\noindent
where $C$ is a complex constant, this has the effect of associating to
each link of the integral-differential chain of $w_{p}(z)$ a complex plane
of constants $C$ in which each point corresponds to a function $w(z)$. In
particular, all constant functions are associated to the null function,
and therefore to a complex plane of constants at each link of the null
chain.

We will now establish a general scheme for the classification of all
possible singularities of inner analytic functions. This can be done for
analytic functions in general, but we will do it here in a way that is
particularly suited for our inner analytic functions.

\begin{definition}\Colon\label{Def05}
  Classification of Singularities: Soft and Hard
\end{definition}

\noindent
Let $z_{1}$ be a point on the unit circle. A singularity of an inner
analytic function $w(z)$ at $z_{1}$ is a {\em soft singularity} if the
limit of $w(z)$ to that point exists and is finite. Otherwise, it is a
{\em hard singularity}.

\vspace{2.6ex}

\noindent
This is a complete classification of all possible singularities because,
given a point of singularity, either the limit of the function to that
point from within the open unit disk exists, or it does not. There is no
third alternative, and therefore every singularity is either soft or hard.
We may now establish the following important property of soft
singularities.

\begin{property}\Colon\label{Prop1SoHaSi}
  A soft singularity of an inner analytic function $w(z)$ at a point
  $z_{1}$ of the unit circle is necessarily an integrable one.
\end{property}

\noindent
In order to prove this first note that, since the singularity at $z_{1}$
is soft, the function $w(z)$ is defined by continuity there, being
therefore continuous at $z_{1}$. Consider now a curve contained within the
open unit disk, that connects to $z_{1}$ along some direction, that has a
finite length, and which is an otherwise arbitrary curve. We have at once
that $w(z)$ is analytic at all points on this curve except $z_{1}$. It
follows that $w(z)$ is continuous, and thus that $|w(z)|$ is continuous,
everywhere on this curve, {\em including} at $z_{1}$. Hence, the limits of
$|w(z)|$ to all points on this curve exist and are finite positive real
numbers.

We now note that this set of finite real numbers must be bounded, because
otherwise there would be a hard singularity of $w(z)$ somewhere within the
open unit disk, where this function is in fact analytic. We conclude
therefore that over the curve the function $w(z)$ is a bounded continuous
function on a finite-length domain, which implies that $w(z)$ is
integrable in that domain. Therefore, we may state that $w(z)$ is
integrable along arbitrary curves reaching the point $z_{1}$ from strictly
within the open unit disk\footnote{Post-publication note: it is important
  to observe here that, given any curve within the unit disk, from any
  internal point, say $z=0$, to $z_{1}$, the value of the integral does
  not depend on the curve. Given any two such curves, one can see this
  considering a small arc around $z_{1}$ connecting the two curves, using
  the Cauchy-Goursat theorem, and considering the limit in which the arc
  becomes infinitesimal.}, which thus establishes this property.

\vspace{2.6ex}

\noindent
We will now prove a couple of important further properties of the
singularity classification, one for soft singularities and one for hard
singularities. For this purpose, let $z_{1}$ be a point on the unit
circle. Let us discuss first a property of soft singularities, which is
related to angular integration.

\begin{property}\Colon\label{Prop2SoHaSi}
  If an inner analytic function $w(z)$ has a soft singularity at $z_{1}$,
  then the angular primitive of $w(z)$ also has a soft singularity at that
  point.
\end{property}

\noindent
In order to prove this we use the fact that a soft singularity is
necessarily an integrable one. We already know that if $w(z)$ has a
singularity at $z_{1}$, then so does its angular primitive
$w^{-1\ldot}(z)$. If we now consider the angular primitive of $w(z)$ at
$z=z_{1}$, we have

\begin{equation}
  w^{-1\ldot}(z_{1})
  =
  -\ii
  \int_{0}^{z_{1}}dz\,
  \frac{w(z)-w(0)}{z},
\end{equation}

\noindent
where the integral can be taken over any simple curve within the open unit
disk. We already know that the integrand is regular at the origin. Since
the singularity of $w(z)$ at $z_{1}$ is soft, that singularity is
integrable along any simple curve within the open unit disk that goes from
$z=0$ to $z=z_{1}$. Therefore, it follows that this integral exists and is
finite, and thus that $w^{-1\ldot}(z_{1})$ exists and is finite. Since the
function $w^{-1\ldot}(z)$ is thus well defined\footnote{Post-publication
  note: it is important to observe here that the fact that the integral
  does not depend on the curve is implicitly being used here.} at $z_{1}$,
as well as analytic around that point, it follows that the singularity of
$w^{-1\ldot}(z)$ at $z_{1}$ is also soft, which thus establishes this
property.

\vspace{2.6ex}

\noindent
Let us discuss now a property of hard singularities, which is related to
angular differentiation.

\begin{property}\Colon\label{Prop3SoHaSi}
  If an inner analytic function $w(z)$ has a hard singularity at $z_{1}$,
  then the angular derivative of $w(z)$ also has a hard singularity at
  that point.
\end{property}

\noindent
If $w(z)$ has a hard singularity at $z_{1}$, then it is not well defined
there, implying that it is not continuous there, and therefore that it is
also not differentiable there. This clearly implies that the angular
derivative $w^{\ldot}(z)$ of $w(z)$, which we already know to also have a
singularity at $z_{1}$, is not well defined there as well. This in turn
implies that the singularity of $w^{\ldot}(z)$ at $z_{1}$ must be a hard
one.

However, the simplest way to prove this property is to note that it
follows from the previous one, that is, from Property~\ref{Prop2SoHaSi}.
We can prove it by reductio ad absurdum, using the fact that, as we have
already shown, the operations of angular differentiation and angular
integration are inverse operations to each other. If we assume that $w(z)$
has a hard singularity at $z_{1}$ and that $w^{\ldot}(z)$ has a soft
singularity at that point, then we have an inner analytic function, namely
$w^{\ldot}(z)$, that has a soft singularity at $z_{1}$, while its angular
primitive, namely $w(z)$, has a hard singularity at that point. However,
according to Property~\ref{Prop2SoHaSi} this is impossible, since angular
integration always takes a soft singularity to another soft singularity.
This establishes, therefore, that this property holds.

\vspace{2.6ex}

\noindent
The use of the operations of angular differentiation and of angular
integration now leads to a refinement of our general classification of
singularities. We will use them to assign to each singularity a {\em
  degree of softness} or a {\em degree of hardness}. Let $w(z)$ be an
inner analytic function and $z_{1}$ a point on the unit circle, and
consider the following two definitions.

\begin{definition}\Colon\label{Def06}
  Classification of Singularities: Gradation of Soft Singularities
\end{definition}

\noindent
Let us assume that $w(z)$ has a soft singularity at $z_{1}$. If an
arbitrarily large number of successive angular differentiations of $w(z)$
always results in a singularity at $z_{1}$ which is still soft, then we
say that the singularity of $w(z)$ at $z_{1}$ is an {\em infinitely soft}
singularity. Otherwise, if $n_{s}$ is the minimum number of angular
differentiations that have to be applied to $w(z)$ in order for the
singularity at $z_{1}$ to become a hard one, then we define $n_{s}$ as the
{\em degree of softness} of the original singularity of $w(z)$ at
$z_{1}$. Therefore, a degree of softness is an integer
$n_{s}\in\{1,2,3,\ldots,\infty\}$.

\begin{definition}\Colon\label{Def07}
  Classification of Singularities: Gradation of Hard Singularities
\end{definition}

\noindent
Let us assume that $w(z)$ has a hard singularity at $z_{1}$. If an
arbitrarily large number of successive angular integrations of $w(z)$
always results in a singularity at $z_{1}$ which is still hard, then we
say that the singularity of $w(z)$ at $z_{1}$ is an {\em infinitely hard}
singularity. Otherwise, if $n_{h}+1$ is the minimum number of angular
integrations that have to be applied to $w(z)$ in order for the
singularity at $z_{1}$ to become a soft one, then we define $n_{h}$ as the
{\em degree of hardness} of the original singularity of $w(z)$ at $z_{1}$.
Therefore, a degree of hardness is an integer
$n_{h}\in\{0,1,2,3,\ldots,\infty\}$.

\vspace{2.6ex}

\noindent
In order to see that this establishes a complete classification of all
possible singularities, let us examine all the possible outcomes when we
apply angular differentiations and angular integrations to inner analytic
functions. We already saw that, if we apply a angular integration to a
soft singularity, then the result is always another soft
singularity. Similarly we saw that, if we apply a angular differentiation
to a hard singularity, then the result is always another hard
singularity. The two remaining alternatives are the application of a
angular integration to a hard singularity, and the application of a
angular differentiation to a soft singularity. In these two cases the
resulting singularity may be either soft or hard, and the remaining
possibilities were dealt with in Definitions~\ref{Def06} and~\ref{Def07}.
Since this applies to all singularities in all integral-differential
chains, it applies to all possible singularities of all inner analytic
functions.

\vspace{2.6ex}

\noindent
In some cases examples of this classification are well known. For
instance, a simple example of an infinitely hard singularity is any
essential singularity. Examples of infinitely soft singularities are
harder to come by, and they are related to integrable real functions which
are infinitely differentiable but not analytic. A simple example of a hard
singularity with degree of hardness $n\geq 1$ is a pole of order $n$.
Examples of soft singularities are the square root, and products of
strictly positive powers with the logarithm.

If a singularity at a given singular point $z_{1}$ on the unit circle is
either infinitely soft or infinitely hard, then the corresponding
integral-differential chain of singularities contains either only soft
singularities or only hard singularities. If the singularity is neither
infinitely soft nor infinitely hard, then at some point along the
corresponding integral-differential chain the character of the singularity
changes, and from that point on the soft or hard character remains
constant at the new value throughout the rest of the integral-differential
chain in that direction. Therefore, in each integral-differential chain
that does not consist of either only soft singularities or only hard
singularities, there is a single transition between two functions on the
chain where the character of the singularity changes.

Let us examine in more detail the important intermediary case in which we
assign to the singularity the degree of hardness zero, which we will also
describe as that of a {\em borderline hard} singularity.

\begin{definition}\Colon\label{Def08}
  Classification of Singularities: Borderline Hard Singularities
\end{definition}

\noindent
Given an inner analytic function $w(z)$ and a point $z_{1}$ on the unit
circle where it has a hard singularity, if a {\em single} angular
integration of $w(z)$ results in a function $w^{-1\ldot}(z)$ which has at
$z_{1}$ a soft singularity, then we say that the original function $w(z)$
has at $z_{1}$ a {\em borderline hard} singularity, that is, a hard
singularity with degree of hardness zero.

\vspace{2.6ex}

\noindent
We establish now the following important property of borderline hard
singularities.

\begin{property}\Colon\label{Prop1GraSin}
  A borderline hard singularity of an inner analytic function $w(z)$ at
  $z_{1}$ must be an integrable one.
\end{property}

\noindent
This is so because the angular integration of $w(z)$ produces an inner
analytic function $w^{-1\ldot}(z)$ which has at $z_{1}$ a soft
singularity, and therefore is well defined at that point. Since the value
of $w^{-1\ldot}(z)$ at $z_{1}$ is given by an integral of $w(z)$ along a
curve reaching that point, that integral must therefore exist and result
in a finite complex number\footnote{Post-publication note: it is important
  to observe here that we may also conclude that the resulting number does
  not depend on the curve, that is, on the direction along which the curve
  connects to $z_{1}$.}. Therefore, the singularity of $w(z)$ at $z_{1}$
must be an integrable one. We may thus conclude that all borderline hard
singularities are integrable ones, which establishes this property.

\vspace{2.6ex}

\noindent
The transition between a borderline hard singularity and a soft
singularity is therefore the single point of transition of the soft or
hard character of the singularities along the corresponding
integral-differential chain. Starting from a borderline hard singularity,
$n_{s}$ angular integrations produce a soft singularity with degree of
softness $n_{s}$, and $n_{h}$ angular differentiations produce a hard
singularity with degree or hardness $n_{h}$. Note that a strictly positive
degree of softness given by $n$ can be identified with a negative degree
of hardness given by $-n$, and vice-versa. A simple example of a
borderline hard singularity is a logarithmic singularity.

Let us end this section with one more important property of hard
singularities.

\begin{property}\Colon\label{Prop2GraSin}
  The borderline hard singularities are the only hard singularities that
  are integrable.
\end{property}

\noindent
If a hard singularity of $w(z)$ at $z_{1}$ has a degree of hardness of one
or larger, then by angular integration it is mapped to another hard
singularity, the hard singularity of $w^{-1\ldot}(z)$ at $z_{1}$. If the
hard singularity of $w(z)$ were integrable, then $w^{-1\ldot}(z)$ would be
well defined at $z_{1}$, and therefore its singularity would be soft
rather than hard. Since we know that the singularity of $w^{-1\ldot}(z)$
is hard, it follows that the singularity of $w(z)$ cannot be integrable.
In other words, all hard singularities with strictly positive degrees of
hardness are necessarily non-integrable singularities, which establishes
this property.

\vspace{2.6ex}

\noindent
Note that the structure of the integral-differential chains establishes
within the space of proper inner analytic functions what may be described
as a structure of discrete fibers, in which the whole space is decomposed
as a set of non-intersecting discrete linear structures. The same is then
true for the corresponding real objects on the unit circle. One can then
reconstruct the space of all inner analytic functions by associating to
each link of each integral-differential chain a complex plane of constants
to be added to the proper inner analytic function of that link, in order
to get all the non-proper inner analytic functions associated to it. In
terms of the corresponding real functions, this corresponds to the
association to each link of a real line of real constants.

\section{Representation of Integrable Real Functions}\label{Sec03}

In this section we will establish the relation between integrable real
functions and inner analytic functions. When we discuss real functions in
this paper, some properties will be globally assumed for these functions.
These are rather weak conditions to be imposed on these functions, that
will be in force throughout this paper. It is to be understood, without
any need for further comment, that these conditions are valid whenever
real functions appear in the arguments. These weak conditions certainly
hold for any integrable real functions that are obtained as restrictions
of corresponding inner analytic functions to the unit circle. The most
basic condition is that the real functions must be measurable in the sense
of Lebesgue, with the usual Lebesgue measure~\cite{RARudin,RARoyden}. In
essence, this is basic infrastructure to allow the functions to be
integrable.

In order to discuss the other global conditions, we must first discuss the
classification of singularities of a real function. The concept of a
singularity itself is the same as that for a complex function, namely a
point where the function is not analytic. The concept of a removable
singularity is well-known for analytic functions in the complex plane.
What we mean by a removable singularity in the case of real functions on
the unit circle is a singular point such that both lateral limits of the
function to that point exist and result in the same real value, but where
the function has been arbitrarily defined to have some other real value.
This is therefore a point were the function can be redefined by
continuity, resulting in a continuous function at that point. The concepts
of soft and hard singularities are carried in a straightforward way from
the case of complex functions, discussed in Section~\ref{Sec02}, to that
of real functions. The only difference is that the concept of the limit of
the function to the point is now taken to be the real one, along the unit
circle.

The second global condition we will impose is that the functions have no
removable singularities. Since they can be easily eliminated, these are
trivial singularities, which we will simply rule out of our discussions in
this paper. Although the presence of even a denumerably infinite set of
such trivial singularities does not significantly affect the results to be
presented here, their elimination does significantly simplify the
arguments to be presented. It is for this reason, that is, for the sake of
simplicity, that we rule out such irrelevant singularities. In addition to
this we will require, as our third an last global condition, that the
number of hard singularities be finite, and hence that they be all
isolated from one another. There will be no limitation on the number of
soft singularities. In terms of the more immediate characteristics of the
real functions, the relevant requirement is that the number of singular
points where a given real function diverges to infinity be finite.

\vspace{2.6ex}

\noindent
In this section we will prove the following theorem.

\begin{theorem}\Colon\label{Theo01}
  Every integrable real function defined on a finite interval can be
  represented by an inner analytic function, and can be recovered almost
  everywhere by means of the limit to the unit circle of the real part of
  that inner analytic function.
\end{theorem}

\noindent
Given an arbitrary real function defined within an arbitrary finite closed
interval, it can always be mapped to a real function within the periodic
interval $[-\pi,\pi]$, by a simple linear change of variables, so it
suffices for our purposes here to examine only the set of real functions
$f(\theta)$ defined in this standard interval. The interval is then mapped
onto the unit circle of the complex $z$ plane. What happens to the values
of the function at the two ends of the interval when one does this is
irrelevant for our purposes here, but for definiteness we may think that
one attributes to the function at the point $z=-1$ the arithmetic average
of the values of the function at the two ends of the periodic interval.
The further requirements to be imposed on these functions are still quite
weak, namely no more than that they be Lebesgue-integrable in the periodic
interval, so that one can attribute to them a set of Fourier
coefficients~\cite{FSchurchill}.

Since for Lebesgue-measurable functions defined within a compact interval
plain integrability and absolute integrability are equivalent
requirements~\cite{RARudin,RARoyden}, we may assume that the functions are
absolutely integrable, without loss of generality. Note that the functions
do not have to be differentiable or even continuous. They may also be
unlimited, possibly diverging to infinity at some singular points, so long
as they are absolutely integrable. This means, of course, that any hard
singularities that they may have at isolated points must be integrable
singularities, which we may thus characterize as {\em borderline hard}
singularities, in a real sense of the term. This means, in turn, that
although the functions may diverge to infinity at isolated points, their
pairs of asymptotic integrals around these points must still exist and be
finite real numbers.

This in turn means that these borderline hard singularities must be
surrounded by open intervals where there are no other borderline hard
singularities, so that the asymptotic integrals around the singular points
can be well defined and finite. It follows that any existing borderline
hard singularities must be isolated from any other borderline hard
singularities. Note that they do not really have to be isolated
singularities in the usual, strict sense of complex analysis, which would
require that they be isolated from {\em all} other singularities. All that
is required is that the borderline hard singularities be isolated from
each other. Hence the requirement that the number of hard singularities be
finite. Note also that one can have any number of soft singularities, even
an infinite number of them. As we pointed out before, in terms of the
properties of the real functions $f(\theta)$, the important requirement is
that the number of singular points on the unit circle where a given real
function diverges to infinity be finite.

\newpage

\begin{proof}\Colon
\end{proof}

\noindent
With all these preliminaries stated, the first thing that we must do here
is, given an arbitrary integrable real function $f(\theta)$ defined within
the periodic interval $[-\pi,\pi]$, to build from it an analytic function
$w(z)$ within the open unit disk of the complex plane. For this purpose we
will use the Fourier coefficients of the given real function. The Fourier
coefficients~\cite{FSchurchill} are defined by

\noindent
\begin{eqnarray}\label{EQFouCoe}
  \alpha_{0}
  & = &
  \frac{1}{\pi}
  \int_{-\pi}^{\pi}d\theta\,
  f(\theta),
  \nonumber\\
  \alpha_{k}
  & = &
  \frac{1}{\pi}
  \int_{-\pi}^{\pi}d\theta\,
  \cos(k\theta)f(\theta),
  \nonumber\\
  \beta_{k}
  & = &
  \frac{1}{\pi}
  \int_{-\pi}^{\pi}d\theta\,
  \sin(k\theta)f(\theta),
\end{eqnarray}

\noindent
where the set of functions
$\left\{\rule{0em}{2ex}1,\cos(k\theta),\sin(k\theta),
  k\in\{1,2,3,\ldots,\infty\}\right\}$, constitutes the Fourier basis of
functions. Since $f(\theta)$ is absolutely integrable, we have that

\begin{equation}
  \int_{-\pi}^{\pi}d\theta\,
  |f(\theta)|
  =
  2\pi M,
\end{equation}

\noindent
where $M$ is a positive and finite real number, namely the average value,
on the periodic interval $[-\pi,\pi]$, of the absolute value of the
function. If we use the triangle inequalities, it follows therefore that
$\alpha_{0}$ exists and that it satisfies

\noindent
\begin{eqnarray}
  |\alpha_{0}|
  & \leq &
  \frac{1}{\pi}
  \int_{-\pi}^{\pi}d\theta\,
  |f(\theta)|
  \nonumber\\
  & = &
  2M,
\end{eqnarray}

\noindent
that is, it is limited by $2M$. Since the elements of the Fourier basis
are all limited smooth functions, and using again the triangle
inequalities, it now follows that all other Fourier coefficients also
exist, and are also all limited by $2M$,

\noindent
\begin{eqnarray}
  |\alpha_{k}|
  & \leq &
  \frac{1}{\pi}
  \int_{-\pi}^{\pi}d\theta\,
  |\cos(k\theta)||f(\theta)|
  \nonumber\\
  & \leq &
  \frac{1}{\pi}
  \int_{-\pi}^{\pi}d\theta\,
  |f(\theta)|
  \nonumber\\
  & = &
  2M,
  \nonumber\\
  |\beta_{k}|
  & \leq &
  \frac{1}{\pi}
  \int_{-\pi}^{\pi}d\theta\,
  |\sin(k\theta)||f(\theta)|
  \nonumber\\
  & \leq &
  \frac{1}{\pi}
  \int_{-\pi}^{\pi}d\theta\,
  |f(\theta)|
  \nonumber\\
  & = &
  2M,
\end{eqnarray}

\noindent
for all $k$, since the absolute values of the sines and cosines are
limited by one. Given that we have the coefficients $\alpha_{k}$ and
$\beta_{k}$, the construction of the corresponding inner analytic function
is now straightforward. We simply define the set of complex coefficients

\noindent
\begin{eqnarray}\label{EQTayCoe}
  c_{0}
  & = &
  \frac{1}{2}\,\alpha_{0},
  \nonumber\\
  c_{k}
  & = &
  \alpha_{k}
  -
  \ii\beta_{k},
\end{eqnarray}

\noindent
for $k\in\{1,2,3,\ldots,\infty\}$. Note that these coefficients are all
limited by $4M$, since, using once more the triangle inequalities, we have

\noindent
\begin{eqnarray}
  |c_{0}|
  & = &
  \frac{1}{2}\,|\alpha_{0}|
  \nonumber\\
  & \leq &
  M,
  \nonumber\\
  |c_{k}|
  & \leq &
  |\alpha_{k}|
  +
  |\beta_{k}|
  \nonumber\\
  & \leq &
  4M.
\end{eqnarray}

\noindent
We now define a complex variable $z$ associated to $\theta$, using an
auxiliary positive real variable $\rho\geq 0$,

\begin{equation}
  z
  =
  \rho\e{\iii\theta},
\end{equation}

\noindent
where $(\rho,\theta)$ are polar coordinates in the complex $z$ plane. We
then construct the following power series around the origin $z=0$,

\begin{equation}
  S(z)
  =
  \sum_{k=0}^{\infty}
  c_{k}z^{k}.
\end{equation}

\noindent
According to the theorems of complex analysis~\cite{CVchurchill}, where
this power series converges in the complex $z$ plane, it converges
absolutely and uniformly to an analytic function $w(z)$. It then follows
that $S(z)$ is in fact the Taylor series of $w(z)$ around $z=0$. We must
now establish that this series converges within the open unit disk,
whatever the values of the Fourier coefficients, given only that they are
all limited by $4M$. In order to do this we will first prove that the
series $S(z)$ is absolutely convergent, that is, we will establish the
convergence of the corresponding series of absolute values

\begin{equation}
  \overline{S}(z)
  =
  \sum_{k=0}^{\infty}
  |c_{k}|\rho^{k}.
\end{equation}

\noindent
Let us now consider the partial sums of this real series, and replace the
absolute values of the coefficients by their common upper bound,

\noindent
\begin{eqnarray}
  \overline{S}_{n}(z)
  & = &
  \sum_{k=0}^{n}
  |c_{k}|\rho^{k}
  \nonumber\\
  & \leq &
  4M
  \sum_{k=0}^{n} 
  \rho^{k},
\end{eqnarray}

\noindent
where $n\in\{0,1,2,3,\ldots,\infty\}$. This is now the sum of a geometric
progression, so that we have

\begin{equation}
  \overline{S}_{n}(z)
  \leq
  4M\,
  \frac{1-\rho^{n+1}}{1-\rho}.
\end{equation}

\noindent
For $\rho<1$ we may now take the $n\to\infty$ limit of the right-hand
side, without violating the inequality, so that we get the sum of a
geometric series,

\begin{equation}
  \overline{S}_{n}(z)
  \leq
  \frac{4M}{1-\rho}.
\end{equation}

\noindent
For $\rho<1$ the right-hand side is now a positive upper bound for all the
partial sums of the series of absolute values. Therefore, since the
sequence $\overline{S}_{n}(z)$ of partial sums is a monotonically
increasing sequence of real numbers which is bounded from above, it now
follows that this real sequence is necessarily a convergent one.

Therefore the series of absolute values $\overline{S}(z)$ is convergent on
the open unit disk $\rho<1$, which in turn implies that the original
series $S(z)$ is absolutely convergent on that same disk. This then
implies that the series $S(z)$ is simply convergent on that same disk.
Since $S(z)$ is a convergent power series, it converges to an analytic
function on the open unit disk, which we may now name $w(z)$. Since this
is an analytic function within the open unit disk, it is an inner analytic
function, the one that corresponds to the real function $f(\theta)$ on the
unit circle,

\begin{equation}
  f(\theta)
  \longrightarrow
  w(z).
\end{equation}

\noindent
The coefficients $c_{k}$ are now recognized as the Taylor coefficients of
the inner analytic function $w(z)$ with respect to the origin. We have
therefore established that from any integrable real function $f(\theta)$
one can define a unique corresponding inner analytic function $w(z)$. This
completes the first part of the proof of Theorem~\ref{Theo01}.

\vspace{2.6ex}

\noindent
Next we must establish that $f(\theta)$ can be recovered as the limit
$\rho\to 1_{(-)}$, from the open unit disk to the unit circle, of the real
part of $w(z)$, so that we can establish the complete correspondence
between the integrable real function and the inner analytic function,

\begin{equation}
  f(\theta)
  \longleftrightarrow
  w(z).
\end{equation}

\begin{proof}\Colon
\end{proof}

\noindent
We start by writing the coefficients $c_{k}$ in terms of $w(z)$ and
discussing their dependence on $\rho$. Since the complex coefficients
$c_{k}$ are the coefficients of the Taylor series of $w(z)$ around $z=0$,
the Cauchy integral formulas of complex analysis, for the function $w(z)$
and its derivatives, written at $z=0$ for the $k^{\rm th}$ derivative of
$w(z)$, tell us that we have

\begin{equation}\label{EQCauFor}
  c_{k}
  =
  \frac{1}{2\pi\ii}
  \oint_{C}dz\,
  \frac{w(z)}{z^{k+1}},
\end{equation}

\noindent
for all $k$, where $C$ is any simple closed curve within the open unit
disk that contains the origin, which we may now take as a circle centered
at $z=0$ with radius $\rho\in(0,1)$. We now note that, since $w(z)$ is
analytic in the open unit disk, so that the explicit singularity at $z=0$
is the only singularity of the integrand on that disk, by the
Cauchy-Goursat theorem the integral is independent of $\rho$ within the
open unit disk, and therefore so are the complex coefficients $c_{k}$.

It thus follows that the coefficients $c_{k}$ are continuous functions of
$\rho$ inside the open unit disk, and therefore that their $\rho\to
1_{(-)}$ limits exist and have those same constant values. Since we have
the relations in Equation~(\ref{EQTayCoe}), the same is true for the
Fourier coefficients $\alpha_{k}$ and $\beta_{k}$. On the other hand, by
construction these are the same coefficients that were obtained from the
real function $f(\theta)$ on the unit circle, and we may thus conclude
that the coefficients $c_{k}$, $\alpha_{k}$ and $\beta_{k}$, for all $k$,
are all constant with $\rho$ and therefore continuous functions of $\rho$
in the whole {\em closed} unit disk. This means that, at least in the case
of the coefficients, the $\rho\to 1_{(-)}$ limit can be taken trivially.

Let us now establish the fact that $f(\theta)$ and the real part
$u(1,\theta)$ of $w(z)$ at $\rho=1$ have exactly the same set of Fourier
coefficients. We consider first the case of the coefficient $\alpha_{0}$.
If we write the Cauchy integral formula in Equation~(\ref{EQCauFor}) for
the case $k=0$ we get

\begin{equation}
  c_{0}
  =
  \frac{1}{2\pi\ii}
  \oint_{C}dz\,
  \frac{w(z)}{z}.
\end{equation}

\noindent
Recalling that $c_{0}=\alpha_{0}/2$ and writing the integral on the circle
of radius $\rho$ using the integration variable $\theta$ we get

\begin{equation}
  \frac{\alpha_{0}}{2}
  =
  \frac{1}{2\pi}
  \int_{-\pi}^{\pi}d\theta\,
  \left[
    u(\rho,\theta)
    +
    \ii
    v(\rho,\theta)
  \right].
\end{equation}

\noindent
Since $\alpha_{0}$ is real, we conclude that the imaginary part in the
right-hand side must be zero, and thus obtain

\begin{equation}
  \alpha_{0}
  =
  \frac{1}{\pi}
  \int_{-\pi}^{\pi}d\theta\,
  u(\rho,\theta),
\end{equation}

\noindent
thus proving that $\alpha_{0}$, which is the $k=0$ Fourier coefficient of
$f(\theta)$, is also the $k=0$ Fourier coefficient of $u(\rho,\theta)$,
for any value of $\rho$, and thus is, in particular, the $k=0$ Fourier
coefficient of $u(1,\theta)$. This is so because, since the $\rho\to
1_{(-)}$ limit of the coefficient $\alpha_{0}$ can be taken, so can the
limit of the integral in the right-hand side. Note that this shows, in
particular, that $u(1,\theta)$ is an integrable real function. In order to
extend the analysis of the coefficients to the case $k>0$ we must first
derive some preliminary relations. Consider therefore the following
integral, on the same circuit $C$ we used in Equation~(\ref{EQCauFor}),

\begin{equation}
  \oint_{C}dz\,
  w(z)z^{k-1}
  =
  0,
\end{equation}

\noindent
with $k>0$. The integral is zero by the Cauchy-Goursat theorem, since for
$k\geq 1$ the integrand is analytic on the whole open unit disk. As before
we write the integral on the circle of radius $\rho$ using the integration
variable $\theta$, to get

\noindent
\begin{eqnarray}
  \int_{-\pi}^{\pi}d\theta\,
  \left[
    u(\rho,\theta)\cos(k\theta)
    -
    v(\rho,\theta)\sin(k\theta)
  \right]
  \hspace{1em}
  &   &
  \nonumber\\
  +
  \ii
  \int_{-\pi}^{\pi}d\theta\,
  \left[
    u(\rho,\theta)\sin(k\theta)
    +
    v(\rho,\theta)\cos(k\theta)
  \right]
  & = &
  0.
\end{eqnarray}

\noindent
We are therefore left with the two identities involving $u(\rho,\theta)$
and $v(\rho,\theta)$,

\noindent
\begin{eqnarray}\label{EQFouCon}
  \int_{-\pi}^{\pi}d\theta\,
  u(\rho,\theta)\cos(k\theta)
  & = &
  \int_{-\pi}^{\pi}d\theta\,
  v(\rho,\theta)\sin(k\theta),
  \nonumber\\
  \int_{-\pi}^{\pi}d\theta\,
  u(\rho,\theta)\sin(k\theta)
  & = &
  -
  \int_{-\pi}^{\pi}d\theta\,
  v(\rho,\theta)\cos(k\theta),
\end{eqnarray}

\noindent
which are valid for all $k>0$ and for all $\rho\in(0,1)$. If we now write
the integrals of the Cauchy integral formulas in Equation~(\ref{EQCauFor})
explicitly as integrals on $\theta$, we get

\noindent
\begin{eqnarray}
  c_{k}
  & = &
  \frac{\rho^{-k}}{2\pi}
  \left\{
    \int_{-\pi}^{\pi}d\theta\,
    \left[
      u(\rho,\theta)\cos(k\theta)
      +
      v(\rho,\theta)\sin(k\theta)
    \right]
  \right.
  \nonumber\\
  &   &
  \hspace{2.4em}
  \left.
    +
    \ii
    \int_{-\pi}^{\pi}d\theta\,
    \left[
      -
      u(\rho,\theta)\sin(k\theta)
      +
      v(\rho,\theta)\cos(k\theta)
    \right]
  \right\}.
\end{eqnarray}

\noindent
Using the identities in Equation~(\ref{EQFouCon}) in order to eliminate
$v(\rho,\theta)$ in favor of $u(\rho,\theta)$ and recalling that
$c_{k}=\alpha_{k}-\ii\beta_{k}$ we get

\begin{equation}
  \alpha_{k}
  -
  \ii
  \beta_{k}
  =
  \frac{\rho^{-k}}{\pi}
  \int_{-\pi}^{\pi}d\theta\,
  \left[
    u(\rho,\theta)\cos(k\theta)
    -
    \ii
    u(\rho,\theta)\sin(k\theta)
  \right],
\end{equation}

\noindent
so that we have the relations for the Fourier coefficients,

\noindent
\begin{eqnarray}
  \alpha_{k}
  & = &
  \frac{\rho^{-k}}{\pi}
  \int_{-\pi}^{\pi}d\theta\,
  u(\rho,\theta)\cos(k\theta),
  \nonumber\\
  \beta_{k}
  & = &
  \frac{\rho^{-k}}{\pi}
  \int_{-\pi}^{\pi}d\theta\,
  u(\rho,\theta)\sin(k\theta).
\end{eqnarray}

\noindent
Since the $\rho\to 1_{(-)}$ limits of the coefficients in the left-hand
sides can be taken, so can the $\rho\to 1_{(-)}$ limits of the integrals
in the right-hand sides. Therefore, taking the limit we have for the
Fourier coefficients,

\noindent
\begin{eqnarray}
  \alpha_{k}
  & = &
  \frac{1}{\pi}
  \int_{-\pi}^{\pi}d\theta\,
  u(1,\theta)\cos(k\theta),
  \nonumber\\
  \beta_{k}
  & = &
  \frac{1}{\pi}
  \int_{-\pi}^{\pi}d\theta\,
  u(1,\theta)\sin(k\theta),
\end{eqnarray}

\noindent
thus completing the proof that the real functions $u(1,\theta)$ and
$f(\theta)$ have exactly the same set of Fourier coefficients. Note, in
passing, that due to the identities in Equation~(\ref{EQFouCon}) these
same coefficients can also be written in terms of $v(1,\theta)$,

\noindent
\begin{eqnarray}
  \alpha_{k}
  & = &
  \frac{1}{\pi}
  \int_{-\pi}^{\pi}d\theta\,
  v(1,\theta)\sin(k\theta),
  \nonumber\\
  \beta_{k}
  & = &
  -\,
  \frac{1}{\pi}
  \int_{-\pi}^{\pi}d\theta\,
  v(1,\theta)\cos(k\theta),
\end{eqnarray}

\noindent
in which the $\cos(k\theta)$ was exchanged for $\sin(k\theta)$ and the
$\sin(k\theta)$ was exchanged for $-\cos(k\theta)$. In fact, this is one
way to state that $u(1,\theta)$ and $v(1,\theta)$ are two mutually
Fourier-conjugate real functions.

Let us now examine the limit $\rho\to 1_{(-)}$ that allows us to recover
from the real part $u(\rho,\theta)$ of the inner analytic function $w(z)$
the original real function $f(\theta)$. We want to establish that we may
state that

\begin{equation}
  f(\theta)
  =
  \lim_{\rho\to 1_{(-)}}
  u(\rho,\theta)
\end{equation}

\noindent
almost everywhere. Let us prove that $u(1,\theta)$ and $f(\theta)$ must
coincide almost everywhere. Simply consider the real function $g(\theta)$
given by

\begin{equation}
  g(\theta)
  =
  u(1,\theta)-f(\theta),
\end{equation}

\noindent
where

\begin{equation}
  u(1,\theta)
  =
  \lim_{\rho\to 1_{(-)}}
  u(\rho,\theta).
\end{equation}

\noindent
Since it is the difference of two integrable real functions, $g(\theta)$
is itself an integrable real function. However, since the expression of
the Fourier coefficients is linear on the functions, and since
$u(1,\theta)$ and $f(\theta)$ have exactly the same set of Fourier
coefficients, it is clear that all the Fourier coefficients of $g(\theta)$
are zero. Therefore, for the integrable real function $g(\theta)$ we have
that $c_{k}=0$ for all $k$, and thus the inner analytic function that
corresponds to $g(\theta)$ is the identically null complex function
$w_{\gamma}(z)\equiv 0$. This is an inner analytic function which is, in
fact, analytic over the whole complex plane, and which, in particular, is
zero over the unit circle, so that we have\footnote{Post-publication note:
  it is important to observe that there is a subtlety here, because we are
  assuming for the identically zero real function $g(\theta)\equiv 0$ the
  result we are trying to prove for all integrable real functions. While
  it is immediately clear that these two identically zero functions do
  correspond to each other, the fact that the identically zero real
  function is the {\em only} real function associated to the identically
  zero inner analytic function is equivalent to the statement of the
  completeness of the Fourier basis. Therefore this proof must remain here
  subject to the fact that the Fourier basis is in fact complete.}
$g(\theta)\equiv 0$. Note, in particular, that the $\rho\to 1_{(-)}$
limits exist at {\em all} points of the unit circle in the case of the
inner analytic function associated to $g(\theta)$. Since our argument is
based on the Fourier coefficients $\alpha_{k}$ and $\beta_{k}$, which in
turn are given by integrals involving these functions, we can conclude
only that

\begin{equation}
  f(\theta)
  =
  \lim_{\rho\to 1_{(-)}}
  u(\rho,\theta)
\end{equation}

\noindent
is valid {\em almost everywhere} over the unit circle. Therefore, we have
concluded that the $\rho\to 1_{(-)}$ limit of $w(z)$ exists and that the
limit of its real part $u(\rho,\theta)$ results in the values of
$f(\theta)$ almost everywhere. This concludes the proof of
Theorem~\ref{Theo01}.

\vspace{2.6ex}

\noindent
Regarding the fact that the proof above is valid only almost everywhere,
it is possible to characterize, up to a certain point, the set of points
where the recovery of the real function $f(\theta)$ may fail, using the
character of the possible singularities of the corresponding inner
analytic function $w(z)$. Wherever $w(z)$ is either analytic or has only
soft singularities on the unit circle, the $\rho\to 1_{(-)}$ limit exists,
and therefore the values of $f(\theta)$ can be recovered. At points on the
unit circle where $f(\theta)$ has hard singularities, $w(z)$ necessarily
also has hard singularities, and therefore the limit does not exist and
thus the values of $f(\theta)$ cannot be recovered. However, in this case
this fact is irrelevant, since $f(\theta)$ is not well defined at these
points to begin with. In any case, by hypothesis there can be at most a
finite number of such points, which therefore form a zero-measure set.

Therefore, the only points where $f(\theta)$ may exist but not be
recoverable from the real part of $w(z)$ are those singular points on the
unit circle where $f(\theta)$ has a {\em soft} singularity, in the real
sense of the term, while $w(z)$ has a {\em hard} singularity, in the
complex sense of the term. In principle this is possible because the
requirement for a singularity to be soft in the complex case is more
restrictive than the corresponding requirement in the real case. For a
singularity to be soft in the real case it suffices that the limits of the
function to the point exist and be the same only along two directions,
coming from either side along the unit circle, but for the singularity to
be soft in the complex case the limits must exist and be the same along
{\em all} directions.

It is indeed possible for a hard complex singularity on the unit circle to
be so oriented that the limit exists along the two particular directions
along the unit circle, but does not exist along other directions. For
example, consider the rather pathological real function

\begin{equation}
  f(\theta)
  =
  \theta\sin\!\left(\frac{\pi^{2}}{\theta}\right),
\end{equation}

\noindent
for $-\pi\leq\theta<0$ and $0<\theta\leq\pi$. It is well known that this
function has an essential singularity at $\theta=0$ in the complex
$\theta$ plane, which is an infinitely hard singularity. However, if
defined by continuity at $\theta=0$ the function is continuous there, and
therefore the singularity at $\theta=0$ is a soft one in the real sense of
the term. The function is also continuous at all other points on the unit
circle. We now observe that, despite having an infinitely hard complex
singularity at $\theta=0$, this is a limited real function on a finite
domain and therefore an integrable real function, which means that we may
still construct an inner analytic function that corresponds to it.
Presumably, this inner analytic function also has an essential singularity
at the point $z=1$, which corresponds to $\theta=0$ on the unit circle.
This fact would then prevent us from obtaining the value of the function
at $\theta=0$ as the $\rho\to 1_{(-)}$ limit of the real part of that
inner analytic function.

The mere fact that one can establish that there is a well defined inner
analytic function for such a pathological real function is in itself
rather unexpected and surprising. Furthermore, one can easily see that
this is not the only example. One can also consider the related example,
this time one in which the singularity is {\em not} soft in the real sense
of the term,

\begin{equation}
  f(\theta)
  =
  \sin\!\left(\frac{\pi^{2}}{\theta}\right),
\end{equation}

\noindent
which is still a limited real function on a finite domain and therefore an
integrable real function, which again means that we may still construct an
inner analytic function that corresponds to it. Many other variations of
these examples can be constructed without too much difficulty.

Excluding all such exceptional cases, we may consider that the recovery of
the real function $f(\theta)$ as the $\rho\to 1_{(-)}$ limit of the real
part of the inner analytic function holds everywhere in the domain of
definition of $f(\theta)$, that is, wherever it is well defined. In order
to exclude all such exceptional cases, all we have to do is to exchange
the condition that there be at most a finite number of hard singularities,
in the real sense of the term, of the integrable real function
$f(\theta)$, for the condition that there be at most a finite number of
hard singularities with finite degrees of hardness, in the complex sense
of the term, of the corresponding inner analytic
function\footnote{Post-publication note: one may consider introducing here
  a further classification of the inner analytic functions. One could say
  that a {\em regular} inner analytic function is one that has, at all its
  singular points on the unit circle, the same status as the corresponding
  real function, regarding the fundamental analytic properties. By
  contrast, an {\em irregular} inner analytic function would be one that
  fails to have the same status as the corresponding real function,
  regarding one or more of the analytic properties, such as those of
  integrability, continuity, and differentiability. Important examples of
  irregular inner analytic functions would be those associated to singular
  distributions, as well as those associated to the examples of singular
  real functions which were mentioned in this section.} $w(z)$.

Once we have the inner analytic function that corresponds to a given
integrable real function, we may consider the integral-differential chain
to which it belongs. There are two particular cases that deserve mention
here. One is that in which the inner analytic functions in the chain do
not have any singularities at all on the unit circle, in which case the
corresponding real functions are all analytic functions of $\theta$ in the
real sense of the term. The other is that in which the inner analytic
functions in the chain have only infinitely soft singularities on the unit
circle, in which case the corresponding real functions are all infinitely
differentiable functions of $\theta$, although they are not analytic. In
this case one can go indefinitely along the chain in either direction
without any change in the soft character of the singularities.

If, on the other hand, one does have borderline hard singularities or soft
singularities with finite degrees of softness, then at some point along
the chain there will be a transition to one or more hard singularities
with strictly positive degrees of hardness, which do not necessarily
correspond to integrable real functions. It can be shown that most of
these singularities are instead associated to either singular
distributions or non-integrable real functions. Their discussion will be
postponed to the aforementioned forthcoming papers.

\section{Behavior Under Analytic Operations}\label{Sub04}

Let us now discuss how the correspondence between inner analytic functions
and integrable real functions behaves under the respective operations of
differentiation and integration, that take us along the corresponding
integral-differential chain. There are two issues here, one being the
existence of the $\rho\to 1_{(-)}$ limit at each point on the unit circle,
the other being whether or not the correspondence between the real
function $f(\theta)$ and the inner analytic function $w(z)$, established
by the construction of the inner analytic function from the integrable
real function, and by the $\rho\to 1_{(-)}$ limit of the real part of the
inner analytic function, survives the operation unscathed.

\vspace{2.6ex}

\noindent
The existence of the limit $\rho\to 1_{(-)}$ hinges on whether the point
at issue is a singular point or not, and then on whether the singularity
at the point is either soft or hard. If a point on the unit circle is {\em
  not} a singularity of the inner analytic function $w(z)$, then the
$\rho\to 1_{(-)}$ limit always exists at that point, no matter how many
angular integrations or angular differentiations are performed on the
inner analytic function, that is, the limit exists throughout the
corresponding integral-differential chain. The same is true if the point
is an infinitely soft singularity of $w(z)$. On the other hand, if it is
an infinitely hard singularity of $w(z)$, then the limit at that point
never exists, in the complex sense, throughout the integral-differential
chain. Note, however, that in some cases the limit may still exist, in the
real sense, along the unit circle.

If the point on the unit circle is a soft singularity of $w(z)$ with a
finite degree of softness $n_{s}$, then the $\rho\to 1_{(-)}$ limit exists
no matter how many angular integrations are performed, since the operation
of angular integration takes soft singularities to other soft
singularities. However, since the operation of angular differentiation may
take soft singularities to hard singularities, the limit will only exist
up to a certain number of angular differentiations, which is given by
$n_{s}-1$. Again, we note that in some cases the limit may still exist
beyond this point, in the real sense, even if it does not exist in the
complex sense.

If the point on the unit circle is a hard singularity of $w(z)$ with a
finite degree of hardness $n_{h}$, including zero, then the $\rho\to
1_{(-)}$ limit does not exist in the complex sense, and will also fail to
exist in that sense for any of the angular derivatives of $w(z)$, since
the operation of angular differentiation takes hard singularities to other
hard singularities. Once more we note that in some cases the limit may
still exist in the real sense, even if it does not exist in the complex
sense. However, since the operation of angular integration may take hard
singularities to soft singularities, the limit will in fact exist after a
certain number of angular integrations of $w(z)$, which is given by
$n_{h}+1$.

Whatever the situation may be, if after a given set of analytic operations
is performed there is at most a finite number of hard singularities, then
the $\rho\to 1_{(-)}$ limit exists almost everywhere, and therefore the
corresponding real function can be recovered at almost all points on the
unit circle. Note, by the way, that the same is true if there is a
denumerably infinite number of hard singularities, so long as they are
{\em not} densely distributed on the unit circle or any part of it, so
that almost all of then can be isolated.

\vspace{2.6ex}

\noindent
The next question is whether or not the relation between the real function
$f(\theta)$ and the inner analytic function $w(z)$ implies the
corresponding relation between the corresponding functions after an
operation of integration or differentiation is applied. This is always
true from a strictly local point of view, since we have shown in
Section~\ref{Sec02} that the operation of angular differentiation on the
open unit disk reduces to the operation of differentiation with respect to
$\theta$ on the unit circle, and that the operation of angular integration
on the open unit disk reduces to the operation of integration with respect
to $\theta$ on the unit circle, up to an integration constant.

There are, however, some global concerns over the unit circle, since the
operations of angular integration and of angular differentiation always
result in {\em proper} inner analytic functions, and there is no
corresponding property of the operations of integration and
differentiation with respect to $\theta$ on the unit circle. Note that the
condition $w(0)=0$, which holds for a proper inner analytic function, is
translated, on the unit circle, to the global condition that the
corresponding real function $f(\theta)$ have zero average value over that
unit circle. This is so because $w(0)=0$ is equivalent to $c_{0}=0$, and
therefore to $\alpha_{0}=0$. However, according to the definition of the
Fourier coefficients in Equation~(\ref{EQFouCoe}), the coefficient
$\alpha_{0}/2$ is equal to that average value.

One way to examine this issue is to use the correspondence between the
Taylor coefficients $c_{k}$ of the inner analytic function $w(z)$ and the
Fourier coefficients $\alpha_{k}$ and $\beta_{k}$ of the integrable real
function $f(\theta)$, which according to our construction of $w(z)$ are
related by the relations in Equation~(\ref{EQTayCoe}). Since we have that

\begin{equation}
  w(z)
  =
  \sum_{k=0}^{\infty}
  c_{k}z^{k},
\end{equation}

\noindent
it follows from the definition of angular differentiation that under that
operation the coefficients $c_{k}$ transform as

\noindent
\begin{eqnarray}
  c_{0}
  & \rightarrow &
  0,
  \nonumber\\
  c_{k}
  & \rightarrow &
  \ii k
  c_{k},
\end{eqnarray}

\noindent
for $k\in\{1,2,3,\ldots,\infty\}$, and it also follows that from the
definition of angular integration that under that operation they transform
as

\noindent
\begin{eqnarray}
  c_{0}
  & \rightarrow &
  0,
  \nonumber\\
  c_{k}
  & \rightarrow &
  -\,
  \frac{\ii}{k}\,
  c_{k},
\end{eqnarray}

\noindent
for $k\in\{1,2,3,\ldots,\infty\}$. If we now look at the Fourier
coefficients, considering their definition in Equation~(\ref{EQFouCoe}),
in the case $k=0$ we have that under differentiation $\alpha_{0}$
transforms as

\noindent
\begin{eqnarray}\label{EQAveVal}
  \alpha_{0}
  & \rightarrow &
  \frac{1}{\pi}
  \int_{-\pi}^{\pi}d\theta\,
  f'(\theta)
  \nonumber\\
  & = &
  \frac{1}{\pi}
  \int_{-\pi}^{\pi}df(\theta),
\end{eqnarray}

\noindent
which is zero so long as $f(\theta)$ is a continuous function, since we
are integrating on a circle. Note that, if $f(\theta)$ is not continuous,
then $f'(\theta)$ is not even a well defined integrable real function, and
we therefore cannot even write the integral, with what we know so far. In
the case $k>0$ we have that under differentiation the Fourier coefficients
transform as

\noindent
\begin{eqnarray}
  \alpha_{k}
  & \rightarrow &
  \frac{1}{\pi}
  \int_{-\pi}^{\pi}d\theta\,
  f'(\theta)\cos(k\theta)
  \nonumber\\
  & = &
  \frac{k}{\pi}
  \int_{-\pi}^{\pi}d\theta\,
  f(\theta)\sin(k\theta)
  \nonumber\\
  & = &
  k\beta_{k},
  \nonumber\\
  \beta_{k}
  & \rightarrow &
  \frac{1}{\pi}
  \int_{-\pi}^{\pi}d\theta\,
  f'(\theta)\sin(k\theta)
  \nonumber\\
  & = &
  -\,
  \frac{k}{\pi}
  \int_{-\pi}^{\pi}d\theta\,
  f(\theta)\cos(k\theta)
  \nonumber\\
  & = &
  -k\alpha_{k},
\end{eqnarray}

\noindent
where we have integrated by parts, noting that the integrated terms are
zero because we are integrating on a circle. We therefore have, so long as
$f(\theta)$ is a continuous function, that

\noindent
\begin{eqnarray}
  \alpha_{0}
  & \rightarrow &
  0,
  \nonumber\\
  \alpha_{k}-\ii\beta_{k}
  & \rightarrow &
  \ii k
  \left(
    \alpha_{k}
    -
    \ii
    \beta_{k}
  \right),
\end{eqnarray}

\noindent
for $k\in\{1,2,3,\ldots,\infty\}$, which are, therefore, the same
transformations undergone by $c_{k}$. In the case of integration
operations the change in $\alpha_{0}$ is indeterminate due to the presence
of an arbitrary integration constant on $\theta$ and, considering once
more the definition of the Fourier coefficients in
Equation~(\ref{EQFouCoe}), we have that for $k>0$ the Fourier coefficients
transform under integration as

\noindent
\begin{eqnarray}
  \alpha_{k}
  & \rightarrow &
  \frac{1}{\pi}
  \int_{-\pi}^{\pi}d\theta\,
  f^{-1\prime}(\theta)\cos(k\theta)
  \nonumber\\
  & = &
  -\,
  \frac{1}{k\pi}
  \int_{-\pi}^{\pi}d\theta\,
  f(\theta)\sin(k\theta)
  \nonumber\\
  & = &
  -\,
  \frac{\beta_{k}}{k},
  \nonumber\\
  \beta_{k}
  & \rightarrow &
  \frac{1}{\pi}
  \int_{-\pi}^{\pi}d\theta\,
  f^{-1\prime}(\theta)\sin(k\theta)
  \nonumber\\
  & = &
  \frac{1}{k\pi}
  \int_{-\pi}^{\pi}d\theta\,
  f(\theta)\cos(k\theta)
  \nonumber\\
  & = &
  \frac{\alpha_{k}}{k},
\end{eqnarray}

\noindent
where we have again integrated by parts, noting once more that the
integrated terms are zero because we are integrating on a circle. We
therefore have, so long as $f(\theta)$ is an integrable function, and so
long as one {\em chooses} the integration constant of the integration on
$\theta$ leading to $f^{-1\prime}(\theta)$ so that $\alpha_{0}$ is mapped
to zero, that

\noindent
\begin{eqnarray}
  \alpha_{0}
  & \rightarrow &
  0,
  \nonumber\\
  \alpha_{k}-\ii\beta_{k}
  & \rightarrow &
  -\,
  \frac{\ii}{k}
  \left(
    \alpha_{k}
    -
    \ii
    \beta_{k}
  \right),
\end{eqnarray}

\noindent
for $k\in\{1,2,3,\ldots,\infty\}$, which are, once more, the same
transformations undergone by $c_{k}$. We therefore see that, from the
point of view of the respective coefficients, the correspondence between
the real function $f(\theta)$ and the inner analytic function $w(z)$
survives the respective analytic operations, so long as the operations
produce integrable real functions on the unit circle, and so long as one
chooses appropriately the integration constant on $\theta$.

Let us discuss the situation in a little more detail, starting with the
operation of integration. As we saw in Property~\ref{Prop1LogInt} of
Section~\ref{Sec02}, angular integration is translated, up to an
integration constant, to integration with respect to $\theta$ on the unit
circle, when we take the $\rho\to 1_{(-)}$ limit. In addition to this,
angular integration never produces new hard singularities out of soft
ones, so that the $\rho\to 1_{(-)}$ limit giving $f^{-1\prime}(\theta)$
exists at all points where those giving $f(\theta)$ exist. We see
therefore that, so long as the integration constant is chosen so as to
satisfy the condition that the function $f^{-1\prime}(\theta)$ have zero
average value over the unit circle, it follows that the correspondence
between the real function $f(\theta)$ and the inner analytic function
$w(z)$ implies the correspondence between the real function
$f^{-1\prime}(\theta)$ and the inner analytic function $w(z)^{-1\ldot}$,

\noindent
\begin{eqnarray}
  f(\theta)
  & \longleftrightarrow &
  w(z)
  \;\;\;\Rightarrow
  \nonumber\\
  f^{-1\prime}(\theta)
  & \longleftrightarrow &
  w^{-1\ldot}(z).
\end{eqnarray}

\noindent
This is valid so long as $f(\theta)$ is an integrable real function. Let
us now discuss the case of the operation of differentiation. As we saw in
Property~\ref{Prop1LogDif} of Section~\ref{Sec02}, angular differentiation
corresponds to differentiation with respect to $\theta$ on the unit
circle, when we take the $\rho\to 1_{(-)}$ limit. However, angular
differentiation can produce new hard singularities out of soft ones, and
can also produce non-integrable hard singularities out of borderline hard
ones. Therefore, we may conclude only that, if all the singularities of
$w(z)$ are soft, which implies that $f(\theta)$ is continuous, then the
correspondence between the real function $f(\theta)$ and the inner
analytic function $w(z)$ does imply the correspondence between the real
function $f'(\theta)$ and the inner analytic function $w^{\ldot}(z)$,

\noindent
\begin{eqnarray}
  f(\theta)
  & \longleftrightarrow &
  w(z)
  \;\;\;\Rightarrow
  \nonumber\\
  f'(\theta)
  & \longleftrightarrow &
  w^{\ldot}(z),
\end{eqnarray}

\noindent
with the exception of the points where $f'(\theta)$ has hard singularities
produced out of soft singularities of $f(\theta)$. Note, however, that
this statement is true even if $w^{\ldot}(z)$ has borderline hard
singularities and therefore $f'(\theta)$ is not continuous.

On the other hand, if $f(\theta)$ is discontinuous at a finite set of
borderline hard singularities of $w(z)$, then $f'(\theta)$ is not even
well defined everywhere, by the usual definition of the derivative of a
real function. In fact, if $w(z)$ has borderline hard singularities then
$w^{\ldot}(z)$ has hard singularities with degrees of hardness equal to
one, which are non-integrable singularities, so that $f'(\theta)$ is not
necessarily an integrable real function. The same is true if the inner
analytic function $w(z)$ has hard singularities with strictly positive
degrees of hardness. The discussion of cases such as these will be given
in the aforementioned forthcoming papers.

Given any inner analytic function that has at most a finite number of
borderline hard singular points and no singularities harder than that, and
the corresponding integral-differential chain, the results obtained here
allow us to travel freely along the integration side of the chain, without
damaging the correspondence between each inner analytic function and the
corresponding real function. The part of the chain where this is valid is
the part to the integration side starting from the link where all the
singularities are either soft or at most borderline hard. What happens
when one travels in the other direction along the chain, starting from
this link, will be discussed in the aforementioned forthcoming papers.

\section{Conclusions and Outlook}\label{Sec05}

We have shown that there is a close and deep relationship between real
functions and complex analytic functions in the unit disk centered at the
origin of the complex plane. This close relation between real functions
and complex analytic functions allows one to use the powerful and
extremely well-known machinery of complex analysis to deal with the real
functions in a very robust way, even if the real functions are very far
from being analytic. For example, the $\rho\to 1_{(-)}$ limit can be used
to define the values of the functions or the values of their derivatives,
at points where these quantities cannot be defined by purely real means.
The concept of inner analytic functions played a central role in the
analysis presented. The integral-differential chains of inner analytic
functions, as well as the classification of singularities of these
functions, which we introduced here, also played a significant role.

One does not usually associate non-differentiable, discontinuous and
unbounded real functions with single analytic functions. Therefore, it may
come as a bit of a surprise that {\em all} integrable real functions are
given by the real parts of certain inner analytic functions on the open
unit disk when one approaches the unit circle. Note, however, that there
are many more inner analytic function within the open unit disk than those
that were examined here, generated by integrable real functions. This
leads to extensions of the relationship between inner analytic functions
and real functions or related objects on the unit circle, which will be
tackled in the aforementioned forthcoming papers.

One important limitation in the arguments presented here is that requiring
that there be only a finite number of borderline hard singularities. It
may be possible, perhaps, to lift this limitation, allowing for a
denumerably infinite set of such integrable singularities. It is probably
not possible, however, to allow for a densely distributed set of such
singularities. Possibly, the limitation that the number of borderline hard
singularities be finite may be exchanged for the limitation that the
number of {\em accumulation points} of a denumerably infinite set of
singular points with borderline hard singularities be finite.

It is quite apparent that the complex-analytic structure presented here
can be used to discuss the Fourier series of real functions, as well as
other aspects of the structure of the Fourier theory of real functions.
The study of the convergence of Fourier series was, in fact, the way in
which this structure was first unveiled. Parts of the arguments that were
presented can be seen to connect to the Fourier theory, such as the role
played by the Fourier coefficients, and the sufficiency of these Fourier
coefficients to represent the functions, which relates to the question of
the completeness of the Fourier basis of functions. This is a rather
extensive discussion, which will be presented in a forthcoming paper.

It is interesting to note that the structure presented here may go some
way towards explaining the rather remarkable fact that physicists usually
operate with singular objects and divergent series in what may seem, from
a mathematical perspective, a rather careless way, while very rarely
getting into serious trouble while doing this. The fact that there is a
robust underlying complex-analytic structure, that in fact explains how
many such murky operations can in fact be rigorously justified, helps one
to understand the unexpected success of this way to operate within the
mathematics used in physics applications. In the parlance of physics, one
may say that the complex-analytic structure within the unit disk functions
as a universal regulator for all real functions, and related singular
objects, which are of interest in physics applications.

We believe that the results presented here establish a new perspective for
the analysis of real functions. The use of the theory of complex analytic
functions makes it a rather deep and powerful point of view. Since complex
analysis and analytic functions constitute such a powerful tool, with so
many applications in almost all areas of mathematics and physics, it is to
be hoped that other applications of the ideas explored here will in due
time present themselves.

\section*{Acknowledgments}

The author would like to thank his friend and colleague Prof. Carlos
Eugênio Imbassay Carneiro, to whom he is deeply indebted for all his
interest and help, as well as his careful reading of the manuscript and
helpful criticism regarding this work.

\bibliography{allrefs_en}\bibliographystyle{ieeetr}

\begin{thebibliography}{1}

\bibitem{CVchurchill}
R.~V. Churchill, {\em Complex Variables and Applications}.
\newblock McGraw-Hill, second~ed., 1960.

\bibitem{FTotCPI}
J.~L. deLyra, ``Fourier theory on the complex plane \mbox{I} -- conjugate pairs
  of fourier series and inner analytic functions,'' {\em arXiv},
  vol.~1409.2582, 2015.

\bibitem{FTotCPII}
J.~L. deLyra, ``Fourier theory on the complex plane \mbox{II} -- weak
  convergence, classification and factorization of singularities,'' {\em
  arXiv}, vol.~1409.4435, 2015.

\bibitem{FTotCPIII}
J.~L. deLyra, ``Fourier theory on the complex plane \mbox{III} -- low-pass
  filters, singularity splitting and infinite-order filters,'' {\em arXiv},
  vol.~1411.6503, 2015.

\bibitem{FTotCPIV}
J.~L. deLyra, ``Fourier theory on the complex plane \mbox{IV} --
  representability of real functions by their fourier coefficients,'' {\em
  arXiv}, vol.~1502.01617, 2015.

\bibitem{FTotCPV}
J.~L. deLyra, ``Fourier theory on the complex plane \mbox{V} --
  arbitrary-parity real functions, singular generalized functions and locally
  non-integrable functions,'' {\em arXiv}, vol.~1505.02300, 2015.

\bibitem{RARudin}
W.~Rudin, {\em Principles of Mathematical Analysis}.
\newblock McGraw-Hill, third~ed., 1976.
\newblock ISBN-13: 978-0070542358, ISBN-10: 007054235X.

\bibitem{RARoyden}
H.~Royden, {\em Real Analysis}.
\newblock Prentice-Hall, third~ed., 1988.
\newblock ISBN-13: 978-0024041517, ISBN-10: 0024041513.

\bibitem{FSchurchill}
R.~V. Churchill, {\em Fourier Series and Boundary Value Problems}.
\newblock McGraw-Hill, second~ed., 1941.

\end{thebibliography}

\end{document}